\documentclass{amsart}
\usepackage{amsbsy,amssymb,amscd,amsmath,amsthm, hyperref,tikz,subcaption, pgf}
\usetikzlibrary{arrows}

\newcommand{\out}[1]{\ }

\newcommand{\NN}{{\mathbb N}}

\newcommand{\RR}{{\mathbb R}}
\newcommand{\CC}{{\mathbb C}}

\renewcommand{\ge}{\geqslant}
\renewcommand{\le}{\leqslant}
\renewcommand{\phi}{\varphi}
\renewcommand{\epsilon}{\varepsilon}

\newcommand{\OUT}[1]{deleted}

\DeclareMathOperator{\isdef}{\overset{def}{=}}


\makeatletter
\newtheorem*{rep@theorem}{\rep@title}
\newcommand{\newreptheorem}[2]{%
\newenvironment{rep#1}[1]{%
 \def\rep@title{#2 \ref{##1}}%
 \begin{rep@theorem}}%
 {\end{rep@theorem}}}
\makeatother

\newtheorem{theorem}{Theorem}[section]
\newreptheorem{theorem}{Theorem}

\theoremstyle{definition}

\theoremstyle{remark}
\newtheorem{remark}[theorem]{Remark}

\numberwithin{equation}{section}

\newtheoremstyle{case}
{3pt}
  {3pt}
  {}
  {}
  {\bfseries}
  {:}
  {.5em}
  {}
\theoremstyle{case}

\numberwithin{subcase}{case}

\begin{document}
\title{A Laplace transform of irregular growth}
\dedicatory{To the memory of Rien Kaashoek}
\author{Jan Wiegerinck}
\address{KdV Institute for Mathematics
\\University of Amsterdam
\\Science Park 105-107
\\P.O. box 94248, 1090 GE Amsterdam
\\The Netherlands}
\email{j.j.o.o.wiegerinck@uva.nl}


\begin{abstract} We give an example of a Laplace transform $\int_\gamma e^{\zeta z} d\mu(\zeta)$ that does not have regular growth. This answers a question in \cite{HL}.

\end{abstract}
\maketitle
\section{Introduction} 

In the present paper I present the proof of the existence of certain entire functions, namely,  Laplace transforms $\int_\gamma e^{zt} d\mu(z)$ of irregular growth, where $\gamma$ is a Jordan \emph{arc}. This was a problem posed by Korevaar in Hayman's original List of Problems, \cite{H}, updated and expanded in \cite{HL}, where it occurs as Problem 7.16 and where it is mentioned that I solved the problem. 

I proved this result as a grad student and stated it in the list of propositions that accompanied my thesis in 1985.  Such propositions are by custom stated without proof. Although the proof is simple, it is unsatisfactory not to have a proof in the literature.

Rien Kaashoek and I met in the late seventies, when I was a student of Jaap Korevaar. Later on Rien became a kind and wise colleague who acted as formal advisor for my PhD student Paul Beneker. In fact, at one  point in the writing process of Beneker's thesis, he gave invaluable advice. My last mathematical contact with Rien was much later when he  asked me to comment on the chapter on entire functions in his joint book with Sjoerd Verduyn Lunel \cite[Chapter 14]{KVL}. This somewhat rekindled my interest in entire functions.

\section{Entire functions of exponential type}
We recall some well known results about entire functions of exponential type. We refer to Levin \cite{L} for the proofs and details. 

An entire function \[f(z)= \sum_{j=0}^\infty\frac{c_j}{j!}z^j\] is called of \emph{exponential type} if for some $C, D>0$ it satisfies $|f(z)|\le Ce^{D|z|}$. From now on $f$ will be an entire function of exponential type. With $M_f(r)=\max_{|z|=r} |f(z)|$, the \emph{type} of $f$ is $\sigma_f=\limsup_{r\to\infty} \log M_f(r)/r$.
The \emph{indicator} of $f$ is the function
\[h_f(\theta)=\limsup_{r\to\infty} \frac{\log|f(re^{i\theta})|}{r},\quad -\pi<\theta\le \pi.
\]
The indicator $h_f$ is continuous and $\sigma_f=\max_\theta h_f(\theta)$.

The \emph{Borel transform} of $f$ is the function 
\[g(z)=\sum_{j=0}^\infty\frac{c_j}{z^{j+1}}.\]
It is holomorphic outside the disc $\{|z|>\sigma_f\}$ (much more is true, cf. \cite{L} Chapter 1, Section 20, but we don't need that). Termwise integration yields
\begin{equation}
    f(z)=\int_\Gamma e^{zs} g(s)\, ds,\label{Borel}
\end{equation}
for any simple closed path of integration $\Gamma$ around $\{|s|>\sigma_f\}$.

\subsection{Zeros} Let $f$ be an entire function of exponential type. The zeros  $\{a_j, j=1,\ldots\}$ of $f$ have the following properties. Set $n_f(r)=\#\{ a_j: |a_j|\le r\}$. Then

\begin{equation}
        \limsup_{r\to \infty}\frac{n_f(r)}{r}\isdef\Delta_f\le e \sigma_f, \label{Z1}
    \end{equation}
\begin{equation}
        \sum_{|a_j|<r} \frac{1}{a_j} \text{is bounded for $r\to \infty$.} \label{Z2}
    \end{equation}

We have 
\begin{theorem}[\cite{L}, Chapter 1, Theorem 15] \label{Hada} If a sequence $(a_j)$ without points of density satisfies \eqref{Z1} and \eqref{Z2}, then the Hadamard canonical product defines an entire function of exponential type.
\end{theorem}
The zero set $\{a_j\}$ of $f$ is said to have an \emph{angular density} if for all but a denumerable set $\theta,\phi$, ($0<\phi-\theta<2\pi$) the limit 
\begin{equation}
    \Delta_f(\theta,\phi)=\lim_{r\to\infty}\frac{n_f(r,\theta,\phi)}{r}
\end{equation}
exists. Here $n_f(r,\theta,\phi)= \#\{a_j: |a_j|\le r, \theta<\arg a_j<\phi\}$.

If $\{a_j\}$ has angular density and $\lim_{r\to \infty}\sum_{|a_j|<r} a_j^{-1}$ exists, we say that the zeros are \emph{regularly distributed}.

\subsection{Regular growth\label{B}} Let $E$ be a Lebesgue measurable subset of $\RR_{>0}$. Recall that $E$ has relative measure 0 if 
\[\lim_{r\to \infty} \frac{\lambda(E\cap (0,r))}{r}=0.\]
Such an $E$ is called an $E^0$ set. We say that an entire function $f$ of exponential type is of \emph{(completely) regular growth} if 
\[ \lim\nolimits_{r\to\infty}^* \frac{\log |f(re^{i\theta})|}{r} \]
exists for all $\theta$. Here $\lim^*$ means that $r$ runs over the complement of an $E^0$ set that is independent of $\theta$.

\begin{theorem}[cf. \cite{L}, Chapter 3, Theorem 4] \label{C} An entire function of exponential type has regular growth if and only if its zeros are regularly distributed.
\end{theorem}

\section{ A Laplace transform of irregular growth}
In this section we prove
\begin{theorem}There exists Laplace transforms of the form 
\begin{equation}
    F(z)=\int_\gamma g(s)e^{zs} ds.
\end{equation}
    Here $\gamma$ is a Jordan arc in $\CC$ and $g$ is holomorphic in a neighborhood of $\gamma$.
\end{theorem}

We start by constructing an entire function of exponential type and irregular growth in the form of a canonical product.

Let $C_k$ be the circle centered at the origin with radius $2^k$, ($k=1,2, \ldots$) and let $a_{k1},\ldots a_{k2^k}$ on $C_k$ be the $2^k$-th roots of unity multiplied by $2^k$. Put $A=\{a_{kj}\}$ and $n(r)=\# \{a\in A, |a|\le r\}$.
We have 
\begin{equation}
\limsup_{r\to\infty}\frac{n(r)}{r}=\lim_{k\to\infty}\frac{n(2^k)}{2^k}=\lim_{k\to\infty}\frac{\sum_{j=1}^k 2^j}{2^k}=2, \label{eq1}
\end{equation}
while
\begin{equation}\liminf_{r\to\infty}\frac{n(r)}{r}=\lim_{k\to\infty}\frac{\sum_{j=1}^{k-1} 2^j}{2^k}=1. \label{eq2}
\end{equation}
Clearly \begin{equation}\label{delta}
    \sum_{|a_{kj}|\le r}\frac{1}{a_{kj}}=0 \quad\text{for all}\quad r\ge 2.
\end{equation}
We also find, with $n(r,\theta, \phi)=\#\{a_j: |a_j|\le r, \theta<\arg a_j<\phi\}$, ($0<\phi-\theta<2\pi$)
\begin{equation}\label{eq3}
    \limsup_{r\to\infty}\frac{n(r,\theta, \phi)}{r}=2\frac{\phi-\theta}{2\pi}\quad  \text{and}\quad  \liminf_{r\to\infty}\frac{n(r,\theta, \phi)}{r}=\frac{\phi-\theta}{2\pi},
\end{equation}
because, with $N(k)$ the number of $k$-th roots of unity in $\{e^{it},  \theta<t<\phi\}$, we have $\lim_{k\to\infty}N(k)/k=\frac{\phi-\theta}{2\pi}$.

We form the canonical product
\begin{equation}\label{canprod}
f(z)=\prod_{a_{nj}\in A}\left(1-\frac z{a_{nj}}\right)e^{z/a_{nj}}=\prod_{n=1}^\infty\left(1-\left(\frac z{2^n}\right)^{2^n}\right).   
\end{equation}

By \eqref{eq1} and \eqref{delta} Theorem \ref{Hada} gives that the canonical product $f$ is a function of exponential type $C$, say. 

Since for every $k\in\NN$ the last product in \eqref{canprod} can be split as  $P_kG_k$ with $P_k$ a polynomial (which does not influence the indicator) and  $G_k$ an entire function with the property that $G_k(z)=G_k(ze^{2\pi i/2^k})$, we infer that the indicator $h(\theta)=C$  for all $\theta$.

By \eqref{eq3} we see that the zeros of $f$ are not regularly distributed. Thus by Theorem \ref{C} the function $f$ has irregular growth. In particular, replacing $f(z)$ by $f(e^{i\theta}z)$ we can assume that $f$ does not grow regularly in the positive real direction. 

Let $g$ be the Borel transform of $f$; $g$ is holomorphic on $\{|z|>C\}$. Let $\gamma$ be the Jordan arc $\gamma(t)=\{(-C-2-t)e^{i\pi t}, -1\le t\le 1\}$ and let $I$ be the interval $I=[-C-3,-C-1]$. Then $\Gamma=\gamma\cup I$ is a simple closed path of integration around $\{|s|<C\}$ and by\eqref{Borel}
\[f(z)=\int_\Gamma g(s)e^{zs} ds.\]

Consider 
\[u(z)= \int_I g(s)e^{zs} ds.\]
Clearly $u$ is bounded in the right half-plane and will tend to 0 along half lines $re^{i\theta}$, $-\pi/2<\theta<\pi/2$, $r>0$. Let $F=f-u$, then
\[F(z)=\int_\gamma g(s)e^{zs} ds\]
does not have regular growth. 
Indeed, if $F$ would have regular growth, then $\lim^*_{r\to\infty} \frac{\log|F(r)|}{r}=C$ and since $\lim_{r\to \infty} u(r) = 0$, also 
$\lim^*_{r\to\infty} \frac{\log|F(r)+u(r)|}{r}=C$, which is not the case.

\begin{remark} In view of Theorem \ref{C} and the demands put on zero sets to be regularly distributed, it is clear that entire functions of exponential type usually do not have regular growth. For such a function the trick with the Borel transform will work, although it may not be obvious which part of the contour $\Gamma$ should be removed.

The product in \eqref{canprod} resembles the example in \cite{L}, Chapter 4, Section 1. A similar computation gives that the type $\sigma_f\le 6$, but we need not use this.
\end{remark}

\section{Funding and/or Competing interests}
No funds, grants, or other support was received. There are no data used as is common in pure mathematics. Results by others that used are properly cited.
The author has no relevant financial or non-financial interests to disclose and there are no competing interests.


\begin{thebibliography}{99}
\bibitem{KVL} M.~A.~Kaashoek \& S.~M.~ Verduyn Lunel, \textit{Completeness Theorems ans Characteristic Matrix Functions.} Birkh\"auser, 2022. xvi+350 pp.

\bibitem{H} W.K. Hayman, \textit{Research problems in Function Theory.}
The Athlone Press [University of London], London, 1967. vii+56 pp.

\bibitem{HL} W.K. Hayman \& E. Lingham, \textit{Research Problems in Function Theory.} Springer, 2019. viii+284 pp.

\bibitem{L} B.Ja. Levin, \textit{Distribution of Zeros of Entire Functions}, Revised edition, Vol 5
Translations of Mathematical Monographs, Amer.~Math.~Soc., 1972. xii+523 pp.
\end{thebibliography}
\end{document}